\overfullrule=0pt
\magnification=1200
\hsize=11.25cm
\vsize=19cm
\hoffset=1cm

\def\sqr#1#2{{\vcenter{\vbox{\hrule height.#2pt \hbox{\vrule width .#2pt height#1pt \kern#1pt \vrule width.#2pt}
\hrule height.#2pt}}}}
\def\square{\mathchoice\sqr64\sqr64\sqr{4.2}3\sqr33} 

\def\S{\smallskip \par}
\def\M{\medskip \par}

\def\BB{\bigskip \bigskip \par}
\def\N{\noindent}

\catcode`@=11
\def\@height{height}
\def\@depth{depth}
\def\@width{width}

\newcount\@tempcnta
\newcount\@tempcntb

\newdimen\@tempdima
\newdimen\@tempdimb

\newbox\@tempboxa

\def\@ifnextchar#1#2#3{\let\@tempe #1\def\@tempa{#2}\def\@tempb{#3}\futurelet
    \@tempc\@ifnch}
\def\@ifnch{\ifx \@tempc \@sptoken \let\@tempd\@xifnch
      \else \ifx \@tempc \@tempe\let\@tempd\@tempa\else\let\@tempd\@tempb\fi
      \fi \@tempd}
\def\@ifstar#1#2{\@ifnextchar *{\def\@tempa*{#1}\@tempa}{#2}}

\def\@whilenoop#1{}
\def\@whilenum#1\do #2{\ifnum #1\relax #2\relax\@iwhilenum{#1\relax 
     #2\relax}\fi}
\def\@iwhilenum#1{\ifnum #1\let\@nextwhile=\@iwhilenum 
         \else\let\@nextwhile=\@whilenoop\fi\@nextwhile{#1}}

\def\@whiledim#1\do #2{\ifdim #1\relax#2\@iwhiledim{#1\relax#2}\fi}
\def\@iwhiledim#1{\ifdim #1\let\@nextwhile=\@iwhiledim 
        \else\let\@nextwhile=\@whilenoop\fi\@nextwhile{#1}}

\newdimen\@wholewidth
\newdimen\@halfwidth
\newdimen\unitlength \unitlength =1pt
\newbox\@picbox
\newdimen\@picht

\def\@nnil{\@nil}
\def\@empty{}
\def\@fornoop#1\@@#2#3{}

\def\@for#1:=#2\do#3{\edef\@fortmp{#2}\ifx\@fortmp\@empty \else
    \expandafter\@forloop#2,\@nil,\@nil\@@#1{#3}\fi}

\def\@forloop#1,#2,#3\@@#4#5{\def#4{#1}\ifx #4\@nnil \else
       #5\def#4{#2}\ifx #4\@nnil \else#5\@iforloop #3\@@#4{#5}\fi\fi}

\def\@iforloop#1,#2\@@#3#4{\def#3{#1}\ifx #3\@nnil 
       \let\@nextwhile=\@fornoop \else
      #4\relax\let\@nextwhile=\@iforloop\fi\@nextwhile#2\@@#3{#4}}

\def\@tfor#1:=#2\do#3{\xdef\@fortmp{#2}\ifx\@fortmp\@empty \else
    \@tforloop#2\@nil\@nil\@@#1{#3}\fi}
\def\@tforloop#1#2\@@#3#4{\def#3{#1}\ifx #3\@nnil 
       \let\@nextwhile=\@fornoop \else
      #4\relax\let\@nextwhile=\@tforloop\fi\@nextwhile#2\@@#3{#4}}

\def\@makepicbox(#1,#2){\leavevmode\@ifnextchar 
   [{\@imakepicbox(#1,#2)}{\@imakepicbox(#1,#2)[]}}

\long\def\@imakepicbox(#1,#2)[#3]#4{\vbox to#2\unitlength
   {\let\mb@b\vss \let\mb@l\hss\let\mb@r\hss
    \let\mb@t\vss
    \@tfor\@tempa :=#3\do{\expandafter\let
        \csname mb@\@tempa\endcsname\relax}%
\mb@t\hbox to #1\unitlength{\mb@l #4\mb@r}\mb@b}}

\def\picture(#1,#2){\@ifnextchar({\@picture(#1,#2)}{\@picture(#1,#2)(0,0)}}

\def\@picture(#1,#2)(#3,#4){\@picht #2\unitlength
\setbox\@picbox\hbox to #1\unitlength\bgroup 
\hskip -#3\unitlength \lower #4\unitlength \hbox\bgroup\ignorespaces}

\def\endpicture{\egroup\hss\egroup\ht\@picbox\@picht
\dp\@picbox\z@\leavevmode\box\@picbox}

\long\def\put(#1,#2)#3{\@killglue\raise#2\unitlength\hbox to \z@{\kern
#1\unitlength #3\hss}\ignorespaces}

\long\def\multiput(#1,#2)(#3,#4)#5#6{\@killglue\@multicnt=#5\relax
\@xdim=#1\unitlength
\@ydim=#2\unitlength
\@whilenum \@multicnt > 0\do
{\raise\@ydim\hbox to \z@{\kern
\@xdim #6\hss}\advance\@multicnt \m@ne\advance\@xdim
#3\unitlength\advance\@ydim #4\unitlength}\ignorespaces}

\def\@killglue{\unskip\@whiledim \lastskip >\z@\do{\unskip}}

\def\thinlines{\let\@linefnt\tenln \let\@circlefnt\tencirc
  \@wholewidth\fontdimen8\tenln \@halfwidth .5\@wholewidth}
\def\thicklines{\let\@linefnt\tenlnw \let\@circlefnt\tencircw
  \@wholewidth\fontdimen8\tenlnw \@halfwidth .5\@wholewidth}

\def\linethickness#1{\@wholewidth #1\relax \@halfwidth .5\@wholewidth}

\def\shortstack{\@ifnextchar[{\@shortstack}{\@shortstack[c]}}

\def\@shortstack[#1]{\leavevmode
\vbox\bgroup\baselineskip-1pt\lineskip 3pt\let\mb@l\hss
\let\mb@r\hss \expandafter\let\csname mb@#1\endcsname\relax
\let\\\@stackcr\@ishortstack}

\def\@ishortstack#1{\halign{\mb@l ##\unskip\mb@r\cr #1\crcr}\egroup}

\def\@stackcr{\@ifstar{\@ixstackcr}{\@ixstackcr}}
\def\@ixstackcr{\@ifnextchar[{\@istackcr}{\cr\ignorespaces}}

\def\@istackcr[#1]{\cr\noalign{\vskip #1}\ignorespaces}

\newif\if@negarg

\def\droite(#1,#2)#3{\@xarg #1\relax \@yarg #2\relax
\@linelen=#3\unitlength
\ifnum\@xarg =0 \@vline 
  \else \ifnum\@yarg =0 \@hline \else \@sline\fi
\fi}

\def\@sline{\ifnum\@xarg< 0 \@negargtrue \@xarg -\@xarg \@yyarg -\@yarg
  \else \@negargfalse \@yyarg \@yarg \fi
\ifnum \@yyarg >0 \@tempcnta\@yyarg \else \@tempcnta -\@yyarg \fi
\ifnum\@tempcnta>6 \@badlinearg\@tempcnta0 \fi
\ifnum\@xarg>6 \@badlinearg\@xarg 1 \fi
\setbox\@linechar\hbox{\@linefnt\@getlinechar(\@xarg,\@yyarg)}%
\ifnum \@yarg >0 \let\@upordown\raise \@clnht\z@
   \else\let\@upordown\lower \@clnht \ht\@linechar\fi
\@clnwd=\wd\@linechar
\if@negarg \hskip -\wd\@linechar \def\@tempa{\hskip -2\wd\@linechar}\else
     \let\@tempa\relax \fi
\@whiledim \@clnwd <\@linelen \do
  {\@upordown\@clnht\copy\@linechar
   \@tempa
   \advance\@clnht \ht\@linechar
   \advance\@clnwd \wd\@linechar}%
\advance\@clnht -\ht\@linechar
\advance\@clnwd -\wd\@linechar
\@tempdima\@linelen\advance\@tempdima -\@clnwd
\@tempdimb\@tempdima\advance\@tempdimb -\wd\@linechar
\if@negarg \hskip -\@tempdimb \else \hskip \@tempdimb \fi
\multiply\@tempdima \@m
\@tempcnta \@tempdima \@tempdima \wd\@linechar \divide\@tempcnta \@tempdima
\@tempdima \ht\@linechar \multiply\@tempdima \@tempcnta
\divide\@tempdima \@m
\advance\@clnht \@tempdima
\ifdim \@linelen <\wd\@linechar
   \hskip \wd\@linechar
  \else\@upordown\@clnht\copy\@linechar\fi}

\def\@hline{\ifnum \@xarg <0 \hskip -\@linelen \fi
\vrule \@height \@halfwidth \@depth \@halfwidth \@width \@linelen
\ifnum \@xarg <0 \hskip -\@linelen \fi}

\def\@getlinechar(#1,#2){\@tempcnta#1\relax\multiply\@tempcnta 8
\advance\@tempcnta -9 \ifnum #2>0 \advance\@tempcnta #2\relax\else
\advance\@tempcnta -#2\relax\advance\@tempcnta 64 \fi
\char\@tempcnta}

\def\vector(#1,#2)#3{\@xarg #1\relax \@yarg #2\relax
\@tempcnta \ifnum\@xarg<0 -\@xarg\else\@xarg\fi
\ifnum\@tempcnta<5\relax
\@linelen=#3\unitlength
\ifnum\@xarg =0 \@vvector 
  \else \ifnum\@yarg =0 \@hvector \else \@svector\fi
\fi
\else\@badlinearg\fi}

\def\@hvector{\@hline\hbox to 0pt{\@linefnt 
\ifnum \@xarg <0 \@getlarrow(1,0)\hss\else
    \hss\@getrarrow(1,0)\fi}}

\def\@vvector{\ifnum \@yarg <0 \@downvector \else \@upvector \fi}

\def\@svector{\@sline
\@tempcnta\@yarg \ifnum\@tempcnta <0 \@tempcnta=-\@tempcnta\fi
\ifnum\@tempcnta <5
  \hskip -\wd\@linechar
  \@upordown\@clnht \hbox{\@linefnt  \if@negarg 
  \@getlarrow(\@xarg,\@yyarg) \else \@getrarrow(\@xarg,\@yyarg) \fi}%
\else\@badlinearg\fi}

\def\@getlarrow(#1,#2){\ifnum #2 =\z@ \@tempcnta='33\else
\@tempcnta=#1\relax\multiply\@tempcnta \sixt@@n \advance\@tempcnta
-9 \@tempcntb=#2\relax\multiply\@tempcntb \tw@
\ifnum \@tempcntb >0 \advance\@tempcnta \@tempcntb\relax
\else\advance\@tempcnta -\@tempcntb\advance\@tempcnta 64
\fi\fi\char\@tempcnta}

\def\@getrarrow(#1,#2){\@tempcntb=#2\relax
\ifnum\@tempcntb < 0 \@tempcntb=-\@tempcntb\relax\fi
\ifcase \@tempcntb\relax \@tempcnta='55 \or 
\ifnum #1<3 \@tempcnta=#1\relax\multiply\@tempcnta
24 \advance\@tempcnta -6 \else \ifnum #1=3 \@tempcnta=49
\else\@tempcnta=58 \fi\fi\or 
\ifnum #1<3 \@tempcnta=#1\relax\multiply\@tempcnta
24 \advance\@tempcnta -3 \else \@tempcnta=51\fi\or 
\@tempcnta=#1\relax\multiply\@tempcnta
\sixt@@n \advance\@tempcnta -\tw@ \else
\@tempcnta=#1\relax\multiply\@tempcnta
\sixt@@n \advance\@tempcnta 7 \fi\ifnum #2<0 \advance\@tempcnta 64 \fi
\char\@tempcnta}

\def\@vline{\ifnum \@yarg <0 \@downline \else \@upline\fi}

\def\@upline{\hbox to \z@{\hskip -\@halfwidth \vrule \@width \@wholewidth
   \@height \@linelen \@depth \z@\hss}}

\def\@downline{\hbox to \z@{\hskip -\@halfwidth \vrule \@width \@wholewidth
   \@height \z@ \@depth \@linelen \hss}}

\def\@upvector{\@upline\setbox\@tempboxa\hbox{\@linefnt\char'66}\raise 
     \@linelen \hbox to\z@{\lower \ht\@tempboxa\box\@tempboxa\hss}}

\def\@downvector{\@downline\lower \@linelen
      \hbox to \z@{\@linefnt\char'77\hss}}

\def\dashbox#1(#2,#3){\leavevmode\hbox to \z@{\baselineskip \z@%
\lineskip \z@%
\@dashdim=#2\unitlength%
\@dashcnt=\@dashdim \advance\@dashcnt 200
\@dashdim=#1\unitlength\divide\@dashcnt \@dashdim
\ifodd\@dashcnt\@dashdim=\z@%
\advance\@dashcnt \@ne \divide\@dashcnt \tw@ 
\else \divide\@dashdim \tw@ \divide\@dashcnt \tw@
\advance\@dashcnt \m@ne
\setbox\@dashbox=\hbox{\vrule \@height \@halfwidth \@depth \@halfwidth
\@width \@dashdim}\put(0,0){\copy\@dashbox}%
\put(0,#3){\copy\@dashbox}%
\put(#2,0){\hskip-\@dashdim\copy\@dashbox}%
\put(#2,#3){\hskip-\@dashdim\box\@dashbox}%
\multiply\@dashdim 3 
\fi
\setbox\@dashbox=\hbox{\vrule \@height \@halfwidth \@depth \@halfwidth
\@width #1\unitlength\hskip #1\unitlength}\@tempcnta=0
\put(0,0){\hskip\@dashdim \@whilenum \@tempcnta <\@dashcnt
\do{\copy\@dashbox\advance\@tempcnta \@ne }}\@tempcnta=0
\put(0,#3){\hskip\@dashdim \@whilenum \@tempcnta <\@dashcnt
\do{\copy\@dashbox\advance\@tempcnta \@ne }}%
\@dashdim=#3\unitlength%
\@dashcnt=\@dashdim \advance\@dashcnt 200
\@dashdim=#1\unitlength\divide\@dashcnt \@dashdim
\ifodd\@dashcnt \@dashdim=\z@%
\advance\@dashcnt \@ne \divide\@dashcnt \tw@
\else
\divide\@dashdim \tw@ \divide\@dashcnt \tw@
\advance\@dashcnt \m@ne
\setbox\@dashbox\hbox{\hskip -\@halfwidth
\vrule \@width \@wholewidth 
\@height \@dashdim}\put(0,0){\copy\@dashbox}%
\put(#2,0){\copy\@dashbox}%
\put(0,#3){\lower\@dashdim\copy\@dashbox}%
\put(#2,#3){\lower\@dashdim\copy\@dashbox}%
\multiply\@dashdim 3
\fi
\setbox\@dashbox\hbox{\vrule \@width \@wholewidth 
\@height #1\unitlength}\@tempcnta0
\put(0,0){\hskip -\@halfwidth \vbox{\@whilenum \@tempcnta < \@dashcnt
\do{\vskip #1\unitlength\copy\@dashbox\advance\@tempcnta \@ne }%
\vskip\@dashdim}}\@tempcnta0
\put(#2,0){\hskip -\@halfwidth \vbox{\@whilenum \@tempcnta< \@dashcnt
\relax\do{\vskip #1\unitlength\copy\@dashbox\advance\@tempcnta \@ne }%
\vskip\@dashdim}}}\@makepicbox(#2,#3)}

\newif\if@ovt 
\newif\if@ovb 
\newif\if@ovl 
\newif\if@ovr 
\newdimen\@ovxx
\newdimen\@ovyy
\newdimen\@ovdx
\newdimen\@ovdy
\newdimen\@ovro
\newdimen\@ovri

\def\@getcirc#1{\@tempdima #1\relax \advance\@tempdima 2pt\relax
  \@tempcnta\@tempdima
  \@tempdima 4pt\relax \divide\@tempcnta\@tempdima
  \ifnum \@tempcnta > 10\relax \@tempcnta 10\relax\fi
  \ifnum \@tempcnta >\z@ \advance\@tempcnta\m@ne
    \else \@warning{Oval too small}\fi
  \multiply\@tempcnta 4\relax
  \setbox \@tempboxa \hbox{\@circlefnt
  \char \@tempcnta}\@tempdima \wd \@tempboxa}

\def\@put#1#2#3{\raise #2\hbox to \z@{\hskip #1#3\hss}}

\def\oval(#1,#2){\@ifnextchar[{\@oval(#1,#2)}{\@oval(#1,#2)[]}}

\def\@oval(#1,#2)[#3]{\begingroup\boxmaxdepth \maxdimen
  \@ovttrue \@ovbtrue \@ovltrue \@ovrtrue
  \@tfor\@tempa :=#3\do{\csname @ov\@tempa false\endcsname}\@ovxx
  #1\unitlength \@ovyy #2\unitlength
  \@tempdimb \ifdim \@ovyy >\@ovxx \@ovxx\else \@ovyy \fi
  \advance \@tempdimb -2pt\relax  %%%% added 7 Dec 89
  \@getcirc \@tempdimb
  \@ovro \ht\@tempboxa \@ovri \dp\@tempboxa
  \@ovdx\@ovxx \advance\@ovdx -\@tempdima \divide\@ovdx \tw@
  \@ovdy\@ovyy \advance\@ovdy -\@tempdima \divide\@ovdy \tw@
  \@circlefnt \setbox\@tempboxa
  \hbox{\if@ovr \@ovvert32\kern -\@tempdima \fi
  \if@ovl \kern \@ovxx \@ovvert01\kern -\@tempdima \kern -\@ovxx \fi
  \if@ovt \@ovhorz \kern -\@ovxx \fi
  \if@ovb \raise \@ovyy \@ovhorz \fi}\advance\@ovdx\@ovro
  \advance\@ovdy\@ovro \ht\@tempboxa\z@ \dp\@tempboxa\z@
  \@put{-\@ovdx}{-\@ovdy}{\box\@tempboxa}%
  \endgroup}

\def\@ovvert#1#2{\vbox to \@ovyy{%
    \if@ovb \@tempcntb \@tempcnta \advance \@tempcntb by #1\relax
      \kern -\@ovro \hbox{\char \@tempcntb}\nointerlineskip
    \else \kern \@ovri \kern \@ovdy \fi
    \leaders\vrule width \@wholewidth\vfil \nointerlineskip
    \if@ovt \@tempcntb \@tempcnta \advance \@tempcntb by #2\relax
      \hbox{\char \@tempcntb}%
    \else \kern \@ovdy \kern \@ovro \fi}}

\def\@ovhorz{\hbox to \@ovxx{\kern \@ovro
    \if@ovr \else \kern \@ovdx \fi
    \leaders \hrule height \@wholewidth \hfil
    \if@ovl \else \kern \@ovdx \fi
    \kern \@ovri}}

\def\circle{\@ifstar{\@dot}{\@circle}}
\def\@circle#1{\begingroup \boxmaxdepth \maxdimen \@tempdimb #1\unitlength
   \ifdim \@tempdimb >15.5pt\relax \@getcirc\@tempdimb
      \@ovro\ht\@tempboxa 
     \setbox\@tempboxa\hbox{\@circlefnt
      \advance\@tempcnta\tw@ \char \@tempcnta
      \advance\@tempcnta\m@ne \char \@tempcnta \kern -2\@tempdima
      \advance\@tempcnta\tw@
      \raise \@tempdima \hbox{\char\@tempcnta}\raise \@tempdima
        \box\@tempboxa}\ht\@tempboxa\z@ \dp\@tempboxa\z@
      \@put{-\@ovro}{-\@ovro}{\box\@tempboxa}%
   \else  \@circ\@tempdimb{96}\fi\endgroup}

\def\@dot#1{\@tempdimb #1\unitlength \@circ\@tempdimb{112}}

\def\@circ#1#2{\@tempdima #1\relax \advance\@tempdima .5pt\relax
   \@tempcnta\@tempdima \@tempdima 1pt\relax
   \divide\@tempcnta\@tempdima 
   \ifnum\@tempcnta > 15\relax \@tempcnta 15\relax \fi    
   \ifnum \@tempcnta >\z@ \advance\@tempcnta\m@ne\fi
   \advance\@tempcnta #2\relax
   \@circlefnt \char\@tempcnta}

%INITIALIZATION
\font\tenln line10
%\font\tencirc circle10
\font\tencirc lcircle10
\font\tenlnw linew10
%\font\tencircw circlew10
\font\tencircw lcirclew10

\thinlines   

\newcount\@xarg
\newcount\@yarg
\newcount\@yyarg
\newcount\@multicnt 
\newdimen\@xdim
\newdimen\@ydim
\newbox\@linechar
\newdimen\@linelen
\newdimen\@clnwd
\newdimen\@clnht
\newdimen\@dashdim
\newbox\@dashbox
\newcount\@dashcnt
\catcode`@=12
%%% Local Variables: 
%%% mode: plain-tex
%%% TeX-master: t
%%% End: 

\catcode`@=11
\font\@linefnt linew10 at 2.4pt
\catcode`@=12

\def\arbreA{\kern-0.4ex
\hbox{\unitlength=.25pt
\picture(60,40)(0,0)
\put(30,0){\droite(0,1){10}}
\put(30,10){\droite(-1,1){30}}
\put(30,10){\droite(1,1){30}}
\endpicture}\kern 0.4ex}

\def\arbreB{\kern-0.4ex
\hbox{\unitlength=.25pt
\picture(60,40)(0,0)
\put(30,0){\droite(0,1){10}}
\put(30,10){\droite(-1,1){30}}
\put(30,10){\droite(1,1){30}}
\put(15,25){\droite(1,1){15}}
\endpicture}\kern 0.4ex}

\def\arbreC{\kern-0.4ex
\hbox{\unitlength=.25pt
\picture(60,40)(0,0)
\put(30,0){\droite(0,1){10}}
\put(30,10){\droite(-1,1){30}}
\put(30,10){\droite(1,1){30}}
\put(45,25){\droite(-1,1){15}}
\endpicture}\kern 0.4ex}

\def\arbreun{\kern-0.4ex
\hbox{\unitlength=.25pt
\picture(60,40)(0,0)
\put(30,0){\droite(0,1){10}}
\put(30,10){\droite(-1,1){30}}
\put(30,10){\droite(1,1){30}}
\put(20,20){\droite(1,1){20}}
\put(10,30){\droite(1,1){10}}
\endpicture}\kern 0.4ex}

\def\arbredeux{\kern-0.4ex
\hbox{\unitlength=.25pt
\picture(60,40)(0,0)
\put(30,0){\droite(0,1){10}}
\put(30,10){\droite(-1,1){30}}
\put(30,10){\droite(1,1){30}}
\put(20,20){\droite(1,1){20}}
\put(30,30){\droite(-1,1){10}}
\endpicture}\kern 0.4ex}

\def\arbretrois{\kern-0.4ex
\hbox{\unitlength=.25pt
\picture(60,40)(0,0)
\put(30,0){\droite(0,1){10}}
\put(30,10){\droite(-1,1){30}}
\put(30,10){\droite(1,1){30}}
\put(50,30){\droite(-1,1){10}}
\put(10,30){\droite(1,1){10}}
\endpicture}\kern 0.4ex}

\def\arbrequatre{\kern-0.4ex
\hbox{\unitlength=.25pt
\picture(60,40)(0,0)
\put(30,0){\droite(0,1){10}}
\put(30,10){\droite(-1,1){30}}
\put(30,10){\droite(1,1){30}}
\put(40,20){\droite(-1,1){20}}
\put(30,30){\droite(1,1){10}}
\endpicture}\kern 0.4ex}

\def\arbrecinq{\kern-0.4ex
\hbox{\unitlength=.25pt
\picture(60,40)(0,0)
\put(30,0){\droite(0,1){10}}
\put(30,10){\droite(-1,1){30}}
\put(30,10){\droite(1,1){30}}
\put(40,20){\droite(-1,1){20}}
\put(50,30){\droite(-1,1){10}}
\endpicture}\kern 0.4ex}

\def\overandunder{\kern-0.4ex
\hbox{\unitlength=.25pt
\picture(440,80)(0,0)
\put(60,0){\droite(0,1){20}}
\put(60,20){\droite(-1,1){60}}
\put(60,20){\droite(1,1){60}}
\put(30,50){\droite(1,1){30}}
\put(20,60){$u$}
\put(50,30){$v$}
\put(-100,0){$u/v=$}

\put(360,0){\droite(0,1){20}}
\put(360,20){\droite(-1,1){60}}
\put(360,20){\droite(1,1){60}}
\put(390,50){\droite(-1,1){30}}
\put(350,30){$u$}
\put(380,60){$v$}
\put(220,0){$u\backslash v=$}

\endpicture}\kern 0.4ex}

\def\hfl#1#2{\smash{\mathop {\hbox to 12mm{\rightarrowfill}} \limits^{\scriptstyle #1}_{\scriptstyle #2}}}

\font\grand=cmr10 at 12pt

\centerline {\grand Order structure on the algebra of permutations}

 \centerline {\grand and of planar binary trees} 
\vskip 2cm

\centerline {\bf Jean-Louis LODAY and Mar\' \i a O. RONCO}
\vskip 2cm

\N {\bf Abstract}.  Let $X_n$ be either the symmetric group on $n$ letters, the set of planar binary $n$-trees
or the set of vertices of the $(n-1)$-dimensional cube. We show that, in each case, the graded 
associative product on $\bigoplus_{n\geq 0} K[X_n]$ can be described explicitly from the weak Bruhat order
 on $X_n$.
\BB

\N {\bf Introduction.} Let $S_n$ be the symmetric group acting on $n$ letters. In [4] Malvenuto
 and Reutenauer showed that the shuffle product induces a graded associative product on the graded 
space $K[S_{\infty }]:= \oplus_{n\geq 0} K[S_n]$ (here $K$ is a field). By using the weak Bruhat order on $S_n$
we  give a close formula for the product of basis elements as follows. Let $\sigma\in S_p$ and $\tau\in S_q$ be
two permutations. We define two operations called respectively `over' and `under':
$$
\sigma / \tau =\sigma \times \tau\in S_{p+q}\quad \hbox {and} \quad \sigma \backslash   \tau 
 =\omega_{p,q}\cdot \sigma \times \tau\in S_{p+q},
$$
where $\omega_{p,q}=(p+1\  p+2\  \cdots p+q\  1\  2\  \cdots p)$.

It turns out that $\sigma / \tau \leq  \sigma \backslash   \tau$ for the weak Bruhat order of $ S_{p+q}$. We prove
 that the product
$*$ on  $K[S_{\infty }]$ is given on the generators by the sum of all permutations in between $\sigma / \tau$ and $
\sigma \backslash   \tau$:
$$
\sigma *\tau = \sum_{\sigma / \tau \leq \omega \leq  \sigma \backslash   \tau}\omega . \leqno (1)
$$

 Let $Y_n$ be the set of planar binary trees with $n$ interior vertices (so the number of elements in $Y_n$
 is the Catalan number
$(2n)! \over n!(n+1)!$). In [5] it is shown that there is a graded associative product
on the graded  space $K[Y_{\infty }]:= \oplus_{n\geq 0} K[Y_n]$ induced by the ``dendriform algebra" structure of
$K[Y_{\infty }]$.
We  give a close formula for the product of basis elements as
follows. There is a partial order on  $Y_n$ induced by $\arbreB < \arbreC$. Let
$u\in Y_p$ and
$v\in Y_q$ be two planar binary trees. We define two operations called respectively `over' and `under' as
follows. The element $u / v \in S_{p+q}$ (resp. $u \backslash v \in S_{p+q}$) is obtained 
by identifying the
root of
$u$ with the left most leaf of $v$ (resp. the right most leaf of $u$ with the root of $v$).
It turns out that $u/ v \leq  u\backslash v$ for the  ordering of $ Y_{p+q}$. We prove that the
product
$*$ on  $K[Y_{\infty }]$ is given on the generators by 
$$
u* v= \sum_{u/ v\leq t \leq  u\backslash v} t . \leqno (2)
$$
In this setting, the operations `over' and `under' are used in [2].
\M

 Let $Q_n= \{\pm 1\}^{n-1}$. There is a graded associative product
on the graded vector space $K[Q_{\infty }]:= \oplus_{n\geq 0} K[Q_n]$ where $ K[Q_n]$ is identified with the Solomon
algebra (cf. [4] and [6]). It is in fact a Hopf algebra whose dual is sometimes called the algebra of
noncommutative symmetric (or quasi-symmetric) functions  (cf. [3]).   We  give a close formula for the
product of basis elements as follows. There is a partial order on 
$Q_n$ induced by
$-1 < +1$. Let $\epsilon \in Q_p$ and
$\delta \in Q_q$. We define two operations called respectively `over' and `under' as
follows: $\epsilon  / \delta := (\epsilon ,-1, \delta) \in Q_{p+q}$ and $\epsilon  \backslash\delta := (\epsilon
,+1,\delta)
\in S_{p+q}$. It is immediate that $\epsilon / \delta \leq \epsilon \backslash \delta$ for the ordering of $Q_{p+q}$. We
prove that the product
$*$ on  $K[Q_{\infty }]$ is given on the generators by 
$$
\epsilon * \delta= \sum_{\epsilon / \delta\leq \alpha  \leq  \epsilon \backslash \delta}  \alpha
= \epsilon / \delta + \epsilon \backslash \delta .\leqno (3)
$$
\M

In [6] we constructed explicit maps 
$$
S_n \buildrel {\psi_n} \over \to  Y_n \buildrel {\phi_n} \over \to Q_n
$$
and we observed that they are in fact restrictions of cellular maps from the cube to the Stasheff polytope
and the permutohedron. Moreover, we showed that, 
after dualization and linear extension, the maps
$$
K[Q_{\infty}] \buildrel {\phi^*} \over \to  K[Y_{\infty}] \buildrel {\psi^*} \over \to K[S_{\infty}]
$$
are injective homomorphisms of  graded  associative algebras. We take advantage of this result to deduce formulas
(2) and (3) from formula (1).
\M
The content of this paper is as follows. In the first part (section 1, 2 and 3) we deal with the partial orders on
$S_n$, $Y_n$ and $Q_n$ respectively, and we show that the maps $\psi_n$ and $\phi_n$ are compatible with the
orders. In the second part (section 4, 5 and 6) we prove formula (1), (2) and (3). In the case of the symmetric
group and in the case of planar binary trees the algebras $K[S_{\infty }]$ and  $K[Y_{\infty }]$ have a more
refined structure: they are dendriform algebras (cf. [5]). We show that in both cases the products $\prec$ and
$\succ$ can also be formulated in terms of the order structure. 

The reader interested only in the symmetric groups,
can read part 1 and then go directly to part 4. The appendix contains a result which is valid for any Coxeter group
and has its own interest.

\vfill
\eject

\N {\bf 1. Weak Bruhat order on the symmetric group $S_n$.} 
\M

Let $(W,S)$ be the Coxeter group $(S_n, \lbrace s_1,\dots ,s_{n-1}\rbrace )$, where $S_n$ is the symmetric group
acting on $\lbrace 1,\dots ,n\rbrace $, and $s_i$ is the transposition of $i$ and $i+1$. We denote by $\cdot$
the group law of $S_n$ and by $1_n$ the unit. In this section we compare the weak Bruhat order on $S_n$ and the
shuffles by applying the result of the Appendix. We
also introduce in 1.9 the operations `over' $/$ and `under'
$\backslash$ from $S_p\times S_q$ to $S_{p+q}$ that are to be used in the Appendix.
\M

 For any permutation $\omega \in S_n$, its length  $l(\omega)$  is the smallest integer $k$ such that $\omega$ can be
written as a product of $k$ generators:  $\omega = s_{i_1}\cdot s_{i_2}\cdot \dots \cdot s_{i_k}$. 

By definition, a permutation $\sigma \in S_n$ has a {\it descent} at $i$, $1\leq i\leq n-1$, if $\sigma (i)>\sigma
(i+1)$. The set of {\it descents } of a permutation $\sigma $ is $Desc(\sigma ):=\lbrace s_i\mid \sigma \ {\rm has \ a\
descent\ at}\ i\rbrace $. 

For any subset  $J\subseteq \lbrace s_1,\dots ,s_{n-1}\rbrace $ the set 
$$X_J^n:= \lbrace \sigma \in
S_n \mid l(\sigma \cdot s_i)>l(\sigma ),\ {\rm for\ all}\ s_i\in J\rbrace $$
 described in the appendix is the set of
all permutations $\sigma \in S_n$ such that $Desc(\sigma )\subseteq \lbrace s_1,\dots ,s_{n-1}\rbrace \setminus J$.
 In order to simplify the notation, we denote the subset $\lbrace s_1,\dots ,s_{p-1},s_{p+1}, \allowbreak \dots
,s_{p+q-1}\rbrace$ of $\lbrace s_1,\dots ,s_{p+q-1}\rbrace $ by $\lbrace s_p\rbrace ^c$.
The set $X_{\lbrace s_p\rbrace ^c }^{p+q}$ is the set of all $(p,q)$-shuffles of $S_{p+q}$,
denoted by $Sh(p,q)$.
There exists a canonical inclusion $\iota : S_p\times S_q\hookrightarrow S_{p+q}$, which maps the generator $s_i$
of $S_p$ to $s_i$ in  $S_{p+q}$, and the generator $s_j$ of $S_q$ to $s_{j+p}$ in $S_{p+q}$. In other words we let a
permutation of $S_p$ act on $\{ 1, \cdots , p\}$ and we let a permutation of $S_q$ act on  $\{ p+1, \cdots , p+q\}$. In the
sequel we identify $S_p\times S_q$ with its image in $S_{p+q}$.

Observe that, for $J=\lbrace s_p\rbrace ^c $, the standard parabolic subgroup $W_{\lbrace s_p\rbrace
^c}$ is precisely $S_p\times S_q$ in $S_{p+q}$. 

Proposition  A.2 of the Appendix takes the following form for the Coxeter group $S_n$: 
\M

\N {\bf Lemma 1.1} {\it Let $p,q\geq 1$.

\N (a) For any $\sigma \in S_{p+q}$ there exist  unique elements $\xi \in Sh(p,q)$ and $\omega \in S_p\times S_q$  such that
$\sigma =\xi \cdot \omega$.

\N (b) For any $\xi \in Sh(p,q)$ and any $\omega \in S_p\times S_q$ the length of $\xi \cdot \omega \in S_{p+q}$
 is the sum: $l(\xi \cdot \omega )=l(\xi )+l(\omega )$. 

\N (c) There exists a longest element in $Sh(p,q)$, denoted $\xi _{p,q}$, and $\xi _{p,q}=(q+1\ q+2\ \cdots  q+p\ 1\
2\ \cdots q)$.}\hfill
$\square$
\M

\N {\bf Definition 1.2} For $n\geq 1$, the {\it weak ordering} (also called {\it weak Bruhat order}) on $S_n$
 is defined as follows: $$\omega \leq \sigma\ {\rm in}\ S_n,\ {\rm if\ there\ exists}\ \tau \in S_n\ {\rm such\ that}\
\sigma=\tau \cdot \omega\ {\rm with}\ l(\sigma )=l(\tau )+l(\omega ).$$ \S

The set of permutations $S_n$, equipped with the weak ordering is a partially ordered set, with minimal element $1_n$,
and maximal element $\omega _n^0:=(n\ n-1\ \dots \ 2\ 1)$ (cf. [1]). 
\M

For $n\geq 1$ and $1\leq i\leq j\leq n-1$, let $c_{i,j}\in S_n$ be the permutation: 
$$c_{i,j}:=s_i\cdot s_{i+1}\cdot \dots \cdot s_j.
$$ 
Given $\omega \in Sh(p,q)$, it is easy to check that, if $\omega \neq 1_{p+q}$,
then there exist integers $l\geq 0$, $1\leq i_1<i_2< \dots <i_{l+1}\leq p+q-1$, and $i_k\leq p+k-1$,
 for $1\leq k\leq l+1$, such
that: 
$$\omega = c_{i_{l+1},p+l}\cdot c_{i_l,p+l-1}\cdot \dots \cdot c_{i_1,p}.
$$ 

Under this notation one has
$$\xi _{p,q}=c_{q,p+q-1}\cdot c_{q-1,p+q-2}\cdot \dots \cdot c_{1,p}.
$$ 
Corollary A.4 of the Appendix and  Lemma 1.1 imply the following result:
\M

\N {\bf Lemma 1.3} {\it Let $p,q\geq 1$ be two integers. The longest element of the set  $Sh(p,q)$ (all $(p,q)$-shuffles)
 is $\xi_{p,q}$. Moreover, one has
$$Sh(p,q)= \lbrace \omega \in S_{p+q} \mid \omega \leq \xi _{p,q}\rbrace .
$$
{ }\hfill $\square$} 
\M

\N {\bf Lemma 1.4} {\it If $\sigma ,\sigma \rq \in S_p$ and $\tau ,\tau \rq \in S_q$ are permutations verifying
 $\sigma \leq \sigma \rq $ and $\tau \leq \tau \rq $, then $\sigma \times \tau \leq \sigma \rq \times \tau \rq $.} 
\M
\N {\it Proof.} The permutations $\sigma \times \tau $ and $\sigma \rq \times \tau \rq $ belong to the subgroup
 $S_p\times S_q$ of $S_{p+q}$.

Since $\sigma \leq \sigma \rq $ and $\tau \leq \tau \rq $, there exist $\delta \in S_p$ and $\epsilon \in S_q$
 such that $\sigma \rq = \delta \cdot \sigma $ and $\tau \rq =\epsilon \cdot \tau $, with $l(\sigma \rq )=
l(\delta ) +
l(\sigma )$ and $l(\tau \rq )=l(\epsilon )+l(\tau )$. 

One has $\sigma \rq \times \tau \rq =\delta \cdot \sigma \times \epsilon \cdot \tau = (\delta \times \epsilon )
\cdot (\sigma \times \tau )$, with $l(\sigma \rq \times \tau \rq ) = l(\sigma \rq )+l(\tau \rq )= l(\delta \times
\epsilon ) + l(\sigma \times \tau )$.\hfill  $\square $
\M

\N {\bf Lemma 1.5} {\it Let $p$ and $q$ be to nonnegative integers, and let $\sigma \in S_p$ and $\tau \in S_q$
 be two 
permutations. If $\omega _1$ and $\omega _2$ are two elements of $Sh(p,q)$ such that $\omega _1<\omega _2$,
 then $$\omega _1\cdot (\sigma \times \tau )<\omega _2\cdot (\sigma \times \tau ).
$$} 

\N {\it Proof.} Apply Lemma 1.1 to $W=S_{p+q}$ and $S=\lbrace s_1,\dots ,s_{p+q-1}\rbrace $. One has that
 $\sigma \times \tau \in S_p\times S_q$ and $\omega _1,\omega _2\in Sh(p,q)$. The result follows
immediately.

{ }\hfill $\square $
\M

\N {\bf Definition 1.6} The {\it grafting } of $\sigma \in S_p$ and $\tau \in S_q$ is the permutation
$\sigma \vee \tau \in S_{p+q+1}$ given by: 
$$(\sigma \vee \tau )(i):=\cases {
\sigma (i)& if $1\leq i\leq p$,\cr
p+q+1&if $i=p+1$,\cr
\tau (i-p-1)+p&if $p+2\leq i\leq p+q+1$. \cr } $$ 
\M

It is easily seen that,
$$\sigma \vee \tau =(\sigma \times \tau \times 1_1)\cdot s_{q+p}\cdot s_{q+p-1}\cdot \dots \cdot s_{p+1},
$$ 
for $\sigma \in S_p$ and $\tau \in S_q$. 
\M

\N {\bf Lemma 1.7} {\it If $\sigma \leq \sigma \rq $ in $S_p$ and $\tau \leq \tau \rq $ in $S_q$,
 then $\sigma \vee \tau \leq \sigma \rq \vee \tau \rq $ in $S_{p+q+1}$.} 
\M

\N {\it Proof.} Suppose $\sigma \rq = \epsilon \cdot \sigma $ and $\tau \rq =\delta \cdot \tau $,
 for some $\epsilon \in S_p$ and $\delta \in S_q$ such that $l(\sigma \rq )=l(\epsilon )+l(\sigma )$
 and $l(\tau \rq
)=l(\delta )+l(\tau )$. Clearly, $\sigma \rq \vee \tau \rq = (\epsilon \times \delta \times 1_1)\cdot 
(\sigma \vee
\tau)$. 

The permutations $\sigma \times \tau \times 1_1$ and $\sigma \rq \times \tau \rq \times 1_1$ belong
 to the subgroup $S_{p+q}\times S_1$ of $S_{p+q+1}$. 

It is immediate to check that $l((s_{p+1}\cdot \dots \cdot s_{p+q})\cdot s_i)>l(s_{p+1}\cdot \dots \cdot s_{p+q})$,
 for
any $1\leq i\leq p+q-1$, that is $s_{p+1}\cdot \dots \cdot s_{p+q}\in Sh(p+q,1)$. Since $s_{p+1} \cdots s_{p+q}   
= c_{p+1, p+q}$, by Lemma 1.1 one has 
$l((s_{p+1}\cdot \dots \cdot s_{p+q})\cdot \omega ) = q+l(\omega )$, for any $\omega \in 
S_{p+q}\times S_1$.
So 
$$l(\sigma \rq \vee \tau \rq)=l((\sigma \rq \vee \tau \rq )^{-1})=l(\sigma \rq )+l(\tau \rq ) + q=
$$ 
$$= l(\epsilon )
+l(\delta )+ l(\sigma )+l(\tau ) +q=l(\epsilon \times \delta \times 1_1)+l(\sigma \vee \tau ).\qquad \square 
$$ 

\N {\bf Proposition 1.8} {\it Let $\sigma \in S_n$ be a permutation such that $\sigma (i)=n$,
 for some $1\leq i\leq n$. There exist unique elements $\sigma^l \in S_{i-1}$, $\sigma^r \in S_{n-i}$ and
 $\gamma \in Sh(i-1, n-i)$ such that: 
$$\sigma = (\gamma \times 1_1)\cdot (\sigma^l \vee \sigma^r ).
$$} 

\N {\it Proof.} Since $\sigma (i)=n$, the element $\sigma $ may be written as $\sigma =\sigma \rq \cdot s_{n-1}
\cdot s_{n-2}\cdot \dots \cdot s_i$, with $\sigma \rq \in S_{n-1}\times S_1$ and $l(\sigma
)=l(\sigma \rq ) +n-i$. 

Lemma 1.1 implies that there exist unique elements $\epsilon \in Sh(i-1,n-i+1)$ and $\delta \in 
W_{\lbrace s_{i-1}\rbrace 
^c}$, such that $\sigma \rq = \epsilon \cdot \delta $. Since the permutation $s_{n-1}$ does not appear 
in a reduced expression of $\sigma \rq $, the following assertions hold:

\N - the element $\epsilon$ is of the form  $\epsilon =\gamma \times 1_1$ for some 
$\gamma \in Sh(i-1,n-i)$.

\N - the element $\delta$ belongs to $  S_{i-1}\times S_{n-i}\times S_1$. So, $\delta =
\sigma^l \times \sigma^r \times 1_1$, for unique $\sigma^l \in S_{i-1}$ and $\sigma^r \in S_{n-i}$. 

Finally, we get that $\sigma =(\gamma \times 1_1)\cdot (\sigma^l \vee \sigma^r )$. The unicity of $\gamma $,
 $\sigma^l $ and $\sigma^r $ follows easily.\hfill $\square $
\M

\N {\bf Definition 1.9} For $p, q\geq 1$, the operations `{\it over}' $/$ and `{\it under}'
$\backslash$ from $S_p\times S_q$ to $S_{p+q}$,
 are defined as follows: 
$$\sigma / \tau :=\sigma \times \tau,\quad {\rm and}\quad \sigma \backslash \tau := \xi _{p,q}
\cdot (\sigma \times \tau ),
$$ for $\sigma \in S_p$ and $\tau \in S_q$. 
\M

Since $\sigma \times \tau \in  S_p\times S_q$, for any $\sigma \in S_p$ and $\tau \in S_q$,
 and $\xi _{p,q}\in Sh(p,q)$, the following relation holds:
$$\sigma /\tau \leq \sigma \backslash \tau .$$ 

\N {\bf Lemma 1.10} {\it The operations $/$ and $\backslash $ are associative.} 
\S

\N {\it Proof.} Let $\sigma \in S_p$, $\tau \in S_q$ and $\delta \in S_r$. It is clear that 
$$(\sigma \times \tau )\times \delta =\sigma \times \tau \times \delta =\sigma \times (\tau \times \delta ).
$$ 

The formula above and the equality
$$\displaylines{
\xi _{p+q,r}\cdot (\xi _{p,q}\times 1_r)=\xi _{p,q+r}\cdot (1_p\times \xi _{q,r})\cr
}$$
imply that the operation $\backslash $ is associative too. \hfill $\square $ 
\BB

\N {\bf 2. Weak ordering on the set of planar binary trees.} 
\M

For $n\geq 1$, let $Y_n$ denote the set of planar binary trees with $n$ vertices:
$$
Y_0 = \{ |\},\  Y_1= \{\  \arbreA \  \},\  Y_2=  \{\  
\arbreB ,\arbreC \  \},\   Y_3= \{\  \arbreun ,\arbredeux
,\arbretrois ,\arbrequatre ,\arbrecinq \  \},
$$
and more generally $Y_n := \bigsqcup _{i+j+1=n } Y_i\times Y_j$. 
The {\it grafting } of a $p$-tree $u$ and a $q$-tree $v$ is the
$(p+q+1)$-tree $u\vee v$ obtained by joining the roots of $u$ and $v$ to a new vertex and create a new root.
 For any tree $t$ there exist unique trees $t^l$ and $t^r$ such that $t=t^l\vee t^r$. 
\S

\N {\bf Definition 2.1} Let $\leq $ be the {\it weak ordering} on $Y_n$ generated transitively by the following 
relations: 

\N a) if $u\leq u\rq \in Y_p$ and $v\leq v\rq \in Y_q$, then $u\vee v\leq u\rq \vee v\rq $ in $Y_{p+q+1}$,

\N b) if $u\in Y_p$, $v\in Y_q$ and $w\in Y_r$, then $(u\vee v)\vee w\leq u\vee (v\vee w)$. 
\S

The pair $(Y_n,\leq )$ is a poset.
\M

\N {\bf Definition 2.2}  The operations `{\it over}' $/$ and `{\it under}'
$\backslash$ from $Y_p\times Y_q$ to $Y_{p+q}$ are defined as follows:

--  $u/v$ is the tree obtained by identifying the root of $u$ with the left most leaf of $v$,

--  $u\backslash v$ is the tree obtained by identifying the right most leaf of $u$ with the root of $v$,
$$
\overandunder
$$
It is immediate to check that $/$ and $\backslash $ are associative. 

Equivalently these operations can be defined  
recursively as follows:

\N - $t/ \vert :=t=: \vert \backslash t $ and  $t\backslash \vert :=t=:\vert / t $ for $t\in Y_n$, 

\N - for $u=u^l\vee u^r$ and $v=v^l\vee v^r$ one has 
$$u/ v := (u/ v^l)\vee v^r,\quad {\rm and }\quad u\backslash v := u^l\vee
(u^r\backslash v).
$$
\N {\bf Lemma 2.3} {\it For any  trees $u\in Y_p$ and $v\in Y_q$ one has
 $$u/v\leq u\backslash v.$$}
\N {\it Proof.} This is an immediate consequence of condition (b) of Definition 2.1. \hfill $\square$
\M

The surjective map $\psi _n:S_n\rightarrow Y_n$ considered in [6] is defined as follows: 

\N - $\psi _1(1_1)=\arbreA \in Y_1$,

\N - the image of a permutation $\sigma \in S_n$ is made of two sequences of integers: the sequence on the left of 
$n$ and the sequence on the right of $n$ in $(\sigma (1),\dots ,\sigma (n))$. These permutations are precisely the
ones apearing in Proposition 1.8. Observe that one of them may be empty. By relabelling the integers in each sequence
so that only consecutive integers (starting with 1) appear, one gets two permutations
$\sigma^l$ and $\sigma^r$. For instance
$(341625)$ gives the two sequences
$(341)$ and $(25)$, which, after relabelling, give $(231)$ and $(12)$. By induction $\psi_n(\sigma)$ is defined as 
$\psi_p(\sigma^l)\vee \psi_q(\sigma^r)$.
\M

\N {\bf Notation.} For $n\geq 1$, let  $S_t$ the subset of $S_n$ such
that a permutation $\sigma \in S_n$ belongs to $S_t$ if and only if $\psi _n(\sigma )=t$:
$$ S_t := \psi_n^{-1} (t) \subset S_n.
$$

This subset admits the following description in terms of shuffles.

\N For $n=1$, one has $S_{\ \arbreA }:=\lbrace 1_1\rbrace = S_1.$ 

\N For $n\geq 2$, let  $t=t^l\vee t^r$, with $t^l\in Y_q$ and $t^r\in Y_p$ and $q+p=n-1$. We have 
$$S_t=\lbrace
(\gamma \times 1_1)\cdot (\sigma \vee \tau ) \mid \gamma \in Sh(p,q), \sigma \in S_{t^l}\ {\rm and}\ \tau \in
S_{t^r}\rbrace ,
$$ 
for $1\leq q\leq n-2$. If $q=0$, then $S_{\vert \vee t^r}= \lbrace \vert \vee \sigma, \hbox { for } \sigma \in
S_{t^r}\rbrace $. And, if $q=n-1$, then $S_{t^l\vee \vert }=\lbrace \sigma \vee \vert, \hbox { for } \sigma \in
S_{t^l}\rbrace $.

For instance when $n=2$, $S_{\ \arbreB }:=\lbrace 1_2\rbrace $ and $S_{\ \arbreC}:=\lbrace s_1\rbrace $. 
 \M

\N {\bf Definition 2.4} For $n\geq 0$, let $Min$ and $Max$ be the maps from $Y_n$ into $S_n$ defined as follows: 

\N For $n=1$, $Min(\ \arbreA ):= 1_1 =:Max(\ \arbreA )$. 

\N For $n=2$, $Min(\ \arbreB ):=1_2=:Max(\ \arbreB )$, and $Min(\ \arbreC ):=s_1=:Max(\ \arbreC )$. 

\N For $n\geq 3$, let $t=t^l\vee t^r$, with $t^l\in Y_q$ and $t^r\in Y_p$ and $p+q=n-1$. The permutations
 $Min(t)$ and $Max(t)$ are defined as follows:

\N If $1\leq q\leq n-2$, then $Min(t):=Min(t^l)\vee Min(t^r)$, and $Max(t):=(\xi _{q,p}\times 1_1)\cdot 
(Max(t^l) \vee Max(t^r))$.

\N If $q=0$, then $Min(t):=\xi _{n-1,1}\cdot (1_1\times Min(t^r))$ and $Max(t):=\xi _{n-1,1}\cdot 
(1_1\times Max(t^r))$. If $q=n-1$, then $Min(t):=Min(t^l)\times 1_1$ and $Max(t):=Max(t^l)\times 1_1$. 
\M

Clearly, $Min(t)$ and $Max(t)$ belong to $S_t$, for any tree $t\in Y_n$. \M

\N {\bf Theorem 2.5}Ê {\it Let $n\geq 1$ and $t\in Y_n$. The following equality holds:
$$S_t=\lbrace \omega \in S_n / Min(t)\leq \omega \leq Max(t)\rbrace .$$ }
\S

\N {\it Proof.} Suppose $t=t^l\vee t^r$, with $t^l\in Y_q$, $t^r\in Y_q$ and $n=q+p+1$.

\N Step 1. Let $\gamma $ and $\gamma \rq $ be elements of $Sh(p,q)$ such that $\gamma \leq \gamma \rq $. Suppose that 
$\sigma\leq \sigma \rq $ in $S_p$ and $\tau\leq \tau \rq$ in $S_q$. 

Lemma 1.7 implies that $\sigma \vee \tau \leq \sigma \rq \vee \tau \rq $. Now, $\gamma \times 1_1$ and $\gamma \rq \times 1_1$ 
belong to $Sh(p,q+1)$, and $\sigma \vee \tau $ and $\sigma \rq \vee \tau \rq $ are elements of $S_p\times S_{q+1}$; 
from Lemma 1.1 one gets,
$$(\gamma \times 1_1)\cdot (\sigma \vee \tau )\leq (\gamma \rq \times 1_1) \cdot (\sigma \rq \vee \tau \rq ).
$$

For any $\gamma \in Sh(p,q)$, Lemma 1.3 states that $1_{p+q}\leq \gamma \leq \xi _{p,q}$. It proves that all 
$\omega \in S_t$ verifies $Min(t)\leq \omega \leq Max(t)$.
\S

\N Step 2. Conversely, let $\omega \in S_n$ be such that $Min(t)\leq \omega \leq Max(t)$.

Since $Min(t)\leq \omega$, there exists $\omega _1\in S_n$ such that $\omega =\omega _1\cdot s_{p+q}\cdot \dots \cdot
s_{p+1}$,  with $l(\omega )=l(\omega _1)+l(s_{p+q}\cdots \dots \cdot s_{p+1})=l(\omega _1) + q$.

\N By Lemma 1.1, there exist unique elements $\omega _2\in Sh(p,q+1)$ and $\omega _3\in S_p\times S_{q+1}$,
such that $\omega _1=\omega _2\cdot \omega _3$, with $l(\omega _1)=l(\omega _2)+l(\omega _3)$.

Since $\omega \leq Max(t)$, there exists $\delta \in S_n$ such that 
$$(\xi _{p,q}\times 1_1)\cdot (Max(t^l)\times Max(t^r)\times 1_1)= \delta \cdot \omega _1,
$$
with $l(\xi _{p,q})+l(Max(t^l))+l(Max(t^r))=l(\delta )+l(\omega _1)$. The permutation $s_{p+q}$ does not appear in a reduced 
expression of $\omega _1$. 

So, $\omega _2\in Sh(p,q+1)$ and $\omega _2(n)=n$, which implies that $\omega _2\leq \xi _{p,q}\times 1_1$.

On the other side, the element $\omega _3\in S_p\times S_{q+1}$ and $s_{p+q}$ does not appear in a reduced 
decomposition of $\omega _3$. So, $\omega _3\in S_p\times S_q\times S_1$. Consequently $\omega _3$ is of the form
$\omega _3=
\sigma _4\times \tau _4\times 1_1$, for unique permutations $\sigma _4\in S_p$ and $\tau _4\in S_q$.
 Moreover, the inequalities
$$(Min(t^l)\times Min(t^r)\times 1_1)\cdot s_{p+q}\cdot \dots \cdot s_{p+1}\leq \omega _2\cdot 
(\sigma _4\times \tau _4\times 1_1)\cdot s_{p+q}\cdot \dots \cdot s_{p+1}\leq 
$$
$$\leq (\xi _{p,q}\times 1_1)\cdot (Max(t^l)\times Max(t^r)\times 1_1)\cdot s_{p+q}\cdot 
\dots \cdot s_{p+1}
$$
imply
$$Min(t^l)\times Min(t^r)\times 1_1\leq \omega _2\cdot (\sigma _4\times \tau _4\times 1_1)\leq (\xi _{p,q}\times 1_1)\cdot 
(Max(t^l)\times Max(t^r)\times 1_1).
$$

Since $1_n\leq \omega _2\leq \xi _{p,q}\times 1_1$ in $Sh(p,q+1)$, by applying Lemma 1.5 we get the following 
identities:
$$Min(t^l)\times Min(t^r)\leq \sigma _4\times \tau _4\leq Max(t^l)\times Max(t^r).
$$

The elements $\sigma _4$ and $\tau _4$ verify that $Min(t^l)\leq \sigma _4\leq Max(t^l)$ and $Min(t^r)\leq \tau _4\leq 
Max(t^r)$.

A recursive argument states that $\sigma _4 \in S_{t^l}$ and $\tau _4\in S_{t^r}$, and the proof is over.
\hfill  $\square $
\M

\N {\bf Corollary 2.6} {\it The weak ordering of $S_n$ induces a partial order $\leq _B$ on $Y_n$.
 This order is compatible with $\psi _n:S_n\rightarrow Y_n$: 
$$\sigma \leq \tau \Rightarrow \psi _n(\sigma )\leq
_B \psi _n (\tau ).$$}

\N {\bf Proposition 2.7} {\it The order $\leq _B$ induced by the weak order on $Y_n$
coincides  with the order $\leq $ of Definition 2.1.}
\S

\N {\it Proof.} We want to see that the order $\leq _B$ verifies conditions (a) and (b) of
Definition 2.1. 

Given $t\in Y_n$ and $w\in Y_m$ recall that, for any $\sigma \in S_t$ and any $\tau \in S_w$, the permutation 
$\sigma \vee \tau $ belongs to $S_{t\vee w}$. Lemma 1.7 implies that $\leq _B$ verifies condition (a). 

Let $t\in Y_n$, $u\in Y_r$ and $w\in Y_m$ be three trees. Suppose that $\sigma \in S_t$, $\delta \in S_u$ and 
$\tau \in S_w$. One has that $(\sigma \vee \delta )\vee \tau $ belongs to $S_{(t\vee u)\vee w}$, while $\sigma \vee
(\delta \vee \tau )$ belongs to $S_{t\vee (u\vee w)}$. To prove condition (b), it suffices to check that $(\sigma \vee
\delta )\vee \tau \leq \sigma \vee (\delta \vee \tau )$ in $S_{n+r+m+2}$. 

Now, an easy calculation shows that:
$$(\sigma \vee \delta )\vee \tau = (\sigma \times \delta \times 1_1\times \tau \times 1_1)\cdot s_{n+r+m+1}\cdot 
\dots \cdot s_{n+r+2}\cdot s_{n+r}\cdot \dots \cdot s_{n+1},
$$ 
and
$$\sigma \vee (\delta \vee \tau )=(\sigma \times \delta \times \tau \times 1_2)\cdot s_{n+r+m}\cdot \dots \cdot 
s_{n+r+1}\cdot s_{n+r+m+1}\cdot \dots \cdot s_{n+1}.
$$
We need to show that $(\sigma \times \delta \times 1_1\times \tau \times 1_1)\cdot s_{n+r+m+1}\cdot \dots \cdot 
s_{n+r+2}$ is smaller than $(\sigma \times \delta \times \tau \times 1_2)\cdot s_{n+r+m}\cdot \dots \cdot
s_{n+r+1}\cdot s_{n+r+m+1}\cdot \dots \cdot s_{n+r+1}$. 
\N We use the relation
$$\displaylines {
s_{n+r+m}\cdot \dots \cdot s_{n+r+1}\cdot s_{n+r+m+1}\cdot \dots \cdot s_{n+r+1}=\hfill\cr
\hfill s_{n+r+m+1}\cdot\dots\cdot  s_{n+r+1}\cdot s_{n+r+m+1}\cdot \dots \cdot s_{n+r+2}.\cr
}$$ 
We have to prove that
$$(\sigma \times \delta \times 1_1\times \tau \times 1_1)\leq (\sigma \times \delta \times \tau \times 1_2) \cdot 
s_{n+r+m+1}\cdot \dots \cdot s_{n+r+1};
$$ 
which is a consequence of the formula:
$$(1_1\times \tau \times 1_1) \leq (\tau \times 1_2)\cdot s_{m+1}\cdot \dots \cdot s_1,\ {\rm for\ any}\ \tau \in 
S_m,\ m\geq 1.\ \leqno (2.6.1) 
$$
To prove $(2.6.1)$ it suffices to check that 
$$l(s_{m+1}\cdot \dots \cdot s_1\cdot (1_1\times \tau \times 1_1))= 
m+1+l(1_1\times \tau \times 1_1),
$$ 
because $s_{m+1}\cdot \dots \cdot s_1$ is in $Sh(1,m+1)$ and $1_1\times \tau 
\times
1_1$ belongs to $S_1\times S_{m+1}$. To end the proof, it suffices to observe  that
$$(\tau \times 1_2)\cdot s_{m+1}\cdot \dots \cdot s_1=s_{m+1}\cdot \dots \cdot s_1\cdot (1_1\times \tau \times 1_1),
\ {\rm for\ any}\ \tau \in S_m,\ m\geq 0.\ \square 
$$

\N {\bf Corollary 2.8} {\it The map $\psi_n: S_n \to Y_n$ is a morphism of posets.} \hfill $\square$
\M

\N {\bf Theorem 2.9} {\it Let $\sigma \in S_p$ and $\tau \in S_q$ be two permutations. 
The following equalities hold: 
$$\psi _{p+q}(\sigma / \tau )= \psi _p (\sigma )/ \psi _q(\tau )
 \quad {\rm and }\quad \psi
_{p+q} (\sigma \backslash \tau )=\psi _p (\sigma )\backslash \psi _q(\tau ).$$} \S

\N {\it Proof.}Ê We prove the first formula, the proof of the second one is similar. 

For $q=1$ the result is obvious.

For $q>1$, suppose that $\tau (r) =q$ for some $1\leq r\leq q$. Proposition 1.8 asserts that there exist 
unique $\gamma \in Sh(r-1,q-r)$, $\tau^l \in S_{r-1}$ and $\tau^r\in S_{q-r}$ such that: 
$$\tau =
 (\gamma \times
1_1)\cdot (\tau^l\vee \tau^r).
$$ 

So one has
$$\sigma \times \tau = (1_p\times \gamma \times 1_1)\cdot (\sigma \times \tau^l\times \tau^r\times 1_1)
\cdot s_{p+q-1}\cdot \dots \cdot s_{p+1}.
$$
The inductive hypothesis states that $\psi _{p+r-1}(\sigma \times \tau^l)=\psi _p(\sigma )/ \psi _{r-1}(\tau^l)$.
 It implies that
$$\sigma \times \tau = (1_p\times \gamma \times 1_1)\cdot ((\sigma \times \tau^l)\vee \tau^r),
$$ 
with $\tau^r\in S_{\psi _{q-r}(\tau^r)}$, $\sigma \times \tau^l\in S_{\psi _p(\sigma )/ \psi
_{r-1}(\tau^l)}$ and 

\N $1_p\times \gamma \in Sh(p+r-1,q-r)$. 

By the definition of
$$S_{\psi _p(\sigma )/ \psi _q (\tau )}=S_{(\psi _p(\sigma )/ \psi _{r-1}(\tau^l))\vee \psi _{q-r}(\tau^r)},
$$ 
we have $\sigma \times \tau \in S_{\psi _p(\sigma )/ \psi _q(\tau )}$. So, since $\sigma /\tau =\sigma \times \tau
$, one has
$$\psi _{p+q}(\sigma / \tau )=\psi _p(\sigma )/ \psi _q(\tau ).
$$ 
{ } \hfill $\square $
\BB

\N {\bf 3. Weak ordering on the set of vertices of the hypercube.}
\M

For $n\geq 2$, let $Q_n:=\lbrace +1,-1\rbrace ^{n-1}$ be the set of
vertices of the hypercube. There is a surjective map
$\phi _n:Y_n\rightarrow Q_n$, which is defined as follows:
$\phi _n(t)=(\epsilon _1,\dots ,\epsilon _{n-1})$, where $\epsilon _i$ is
$-1$ when the $i$th leaf of $t$ is
right oriented (more precisely SW-NE), and $+1$ when it is left oriented
(more precisely SE-NW). We take into account only the
interior leaves of
$t$, since the orientation of  the two extreme ones does not depend on $t$.
For instance $\psi (\ \arbreC) = (+1)$ and $\psi (\
\arbreB) = (-1)$. By convention $Q_1= \lbrace (-1)_1\rbrace $ and 
$\psi (\ \arbreA) = (-1)_1$.
\S

We consider $Q_2$ as the partially ordered set $Q_2:=\lbrace -1<+1\rbrace $.
\M

\N {\bf Definition 3.1} The set $Q_n$ of vertices of the hypercube is a
partially ordered set for the order:
$$\epsilon \leq \eta \ {\rm if\ and\ only\ if}\ \epsilon _i\leq \eta _i,\
{\rm for\ all}\ 1\leq i\leq n-1.$$
\M

We denote by $(-1)_n$ the minimal element of $Q_n$, and by $(+1)_n$ its maximal element.
\M

\N {\bf Definition 3.2}  Given an element $\epsilon =(\epsilon _1,\dots
,\epsilon _{p-1})\in Q_p$ and an element
$\eta =(\eta _1,\dots ,\eta _{q-1})\in Q_q$ the {\it grafting} of $\epsilon
$ and $\eta $, denoted $\epsilon \vee \eta $,
is the element of $Q_{p+q+1}$ given by:
$$\epsilon \vee \eta := (\epsilon _1,\dots ,\epsilon _{p-1},-1,+1,\eta
_1,\dots ,\eta _{q-1}).
$$
The operations over $/$ and under $\backslash $ from $Q_p\times Q_q$ to
$Q_{p+q}$  are defined by
$$\displaylines{
\epsilon /\eta := (\epsilon _1,\dots ,\epsilon _{n-1},-1,\eta _1,\dots
,\eta _{m-1}),\cr
 \epsilon \backslash \eta :=(\epsilon _1,\dots ,\epsilon _{n-1},+1,\eta
_1,\dots ,\eta _{m-1}).\cr
}$$

\N {\bf Remark 3.3} It is easily seen that the maps $\phi $ preserve the
operations grafting $\vee $, over $/$, and under $\backslash$. 
\M

\N {\bf Lemma 3.4} {\it Let $t$ be an element of $Y_n$ such that its $i$th
leaf points to the right, for some $2\leq i\leq n-1$.
If $w$ is another tree in $Y_n$ such that $w\leq t$, then the $i$th leaf of
$w$ is right oriented too.}
\S

\N {\it Proof.} The result is obvious for $n\leq 2$.

Since the order $\leq $ on $Y_n$ is  transitively generated by the
relations given in Definition 2.1, it suffices to show
that the assertion is true for the situations described in (a) and (b) of
this Definition.

\N For (a): If $w=w^l\vee w^r$ and $t=t^l\vee t^r$, with $w^l\leq t^l$ and
$w^r\leq t^r$, then the results is
 an immediate consequence of the
inductive hypothesis.

\N For (b): Suppose $w=(u\vee v)\vee s$ and $t=u\vee (v\vee s)$, for some
$u\in Y_p$, $v\in Y_q$ and $s\in Y_r$. If
$q\geq 1$, then the $k$th leaf of $w$ is oriented in the same direction
that the $k$th leaf of $t$, for all
$1\leq k\leq n+1$.

\N If $q=0$, then the $k$th leaf of $w$ is oriented in the same direction
that the $k$th leaf of $t$, for all
$k\neq p+2$. And the $(p+2)$th leaf of $w$ is right oriented, while
$(p+2)$th leaf of $t$ is left oriented. $\square $
\M

\N {\bf Proposition 3.5} {\it For all $n\geq 1$ and all $\epsilon \in Q_n$ there exist two trees in $Y_n$, denoted 
$min(\epsilon )$ and $max(\epsilon )$ respectively, such that the inverse image of $\epsilon $ by
$\phi _n:Y_n\rightarrow Q_n$ satisfies:
$$\phi ^{-1}(\epsilon )=\lbrace t\in Y_n \mid min(\epsilon )\leq t\leq
max(\epsilon )\rbrace .
$$}

\N {\it Proof.} Step 1. The  inverse image $\phi_n^{-1}((-1)_n)$ of the minimal element of $Q_n$
 is the minimal tree $a_n$  of $Y_n$ which has all its leaves pointing to the right. Similarly,
the inverse image $\phi _n^{-1}((+1) _n)$ of the maximal element of $Q_n$ is the maximal tree $z_n$ of $Y_n$ which
has all its leaves   pointing to the left. So, the theorem is obviously true for $\epsilon \in \lbrace (-1)_n,(+1)_n\rbrace $ if
we define:
$$min ((-1)_n):=a_n=:max((-1)_n),\ {\rm and\ }min ((+1)_n):=z_n=:max((+1)_n).$$

\N If $\epsilon \notin \lbrace (-1)_n;(+1)_n\rbrace $, we define $max$ and $min$ recursively, as follows:
\S
\N (a) If $\epsilon _1 =-1$ there exist $k\geq 1$ and $\epsilon \rq \in Q_{n-k}$ such that 
$\epsilon =(-1)_k/ \epsilon \rq $. Define $min(\epsilon ):=a_k/min(\epsilon \rq )$.

\N If $\epsilon_1 =+1$, there exist $k\geq 2$ and $\epsilon \rq \in Q_{n-k}$ such that 
$\epsilon =(+1)_k/\epsilon \rq $. Define $min(\epsilon ):=z_k/min(\epsilon \rq )$.
\S
\N (b) If $\epsilon _{n-1}=-1$, there exist $k\geq 2$ and $\epsilon \rq \in Q_{n-k}$ such that 
$\epsilon =\epsilon \rq \backslash (-1)_k$. Define $max (\epsilon ):= max (\epsilon \rq )\backslash a_k$.

\N If $\epsilon _{n-1}=+1$, there exist $k\geq 1$ and $\epsilon \rq \in Q_{n-k}$ such that 
$\epsilon =\epsilon \rq \backslash (+1)_k$. Define $max (\epsilon ):= max (\epsilon \rq )\backslash z_k$.
\S

\N Step 2. It is easy to prove, by induction on $n$, that if $t\in
\phi _n^{-1}(\epsilon )$, then
$min(\epsilon )\leq t\leq max(\epsilon )$.

Conversely, let $t$ be a tree such that $min(\epsilon )\leq t\leq
max(\epsilon )$. Since $min(\epsilon )\leq t$, Lemma 3.4
implies that the $i$th leaf of $t$ is left oriented, for all $i$ such that
$\epsilon _i=+1$. Similarly,
$t\leq max(\epsilon )$ and Lemma 3.4 imply that the $i$th leaf of $t$ is
right oriented, for all $i$ such that $\epsilon _i=-1$.
So, $t$ belongs to $\phi _n^{-1}(\epsilon )$. \hfill $\square $
\M

\N {\bf Corollary 3.6} {\it For $n\geq 2$, the order of $Y_n$ induces a
partial order $\leq _B$ on $Q_n$.
 This order is compatible with $\phi _n:Y_n\rightarrow Q_n$.}
\M

\N {\bf Proposition 3.7} {\it The order $\leq _B$ of $Q_n$ coincides with
the order $\leq $ of Definition 3.1.}
\S

\N {\it Proof.} If $w$ and $t$ are two trees in $Y_n$ such that $w\leq t$,
then Lemma 3.4 implies that $\phi _n(w)\leq
\phi _n(t)$. It proves that if $\epsilon \leq _B \eta $ in $Q_n$, then
$\epsilon \leq \eta $.

\N To prove that $\epsilon \leq \eta $ in $Q_n$ implies that $\epsilon \leq
_B \eta $, it suffices to show that
$$(\epsilon _1,\dots ,\epsilon _{p-1},-1,\epsilon _{p+1},\dots ,\epsilon
_{n-1})\leq _B ((\epsilon _1,\dots ,\epsilon _p,
+1,\epsilon _{p+2},\dots ,\epsilon _{n-1}),
$$
for all $1\leq p\leq n-1$ and all elements $\epsilon _i\in \lbrace
-1,+1\rbrace $, $1\leq i\leq n-1,\ i\neq p+1$.
Consider the element $\kappa :=(\epsilon _1,\dots ,\epsilon _{p-1})$ in
$Q_p$, and the element
$\rho  :=(\epsilon _{p+1},\dots ,\epsilon _{n-1})\in Q_{n-p}$. Let $t\in
Y_p$ be a tree in $\phi _p^{-1}(\kappa )$ and $w\in
Y_{n-p}$ be a tree in $\phi_{n-p}^{-1}(\rho )$. It is easy to check that
$$\matrix {
\phi _n(t/w)&=&(\epsilon _1,\dots ,\epsilon _{p-1},-1,\epsilon _{p+1},\dots
,\epsilon _{n-1})&\ {\rm and}\cr
\phi _n(t\backslash w)&=&(\epsilon _1,\dots ,\epsilon _{p-1},+1,\epsilon
_{p+1},\dots ,\epsilon _{n-1}).&\hfill \cr }
$$
Since Lemma 2.3 states that $t/w\leq t\backslash w$ in $Y_n$, one gets the
result. $\square $
\M

\N {\bf Corollary 3.8}  {\it The map $\phi _n:Y_n\rightarrow Q_n$ is a
morphism of posets.}\hfill $\square$
\BB

\N {\bf 4. The graded algebra of permutations {\bf Q}$[S_{\infty }]$.}
\M

Consider the graded vector space {\bf Q}$[S_{\infty }]:=\oplus _{n\geq 0}${\bf Q}$[S_n]$, equipped with the shuffle
product $*$ defined by:
$$\sigma *\tau :=\sum _{x\in Sh(p,q)}x\cdot (\sigma \times \tau ),\ {\rm for}\ \sigma \in S_p
\ {\rm and}\ \tau \in S_q.
$$
In [4], C. Malvenuto and C. Reutenauer prove that ({\bf Q}$[S_{\infty }],*)$ is an associative algebra
 over {\bf Q}. We denote by $\overline {{\bf Q}[S_{\infty }]}$ the augmentation ideal.
\M

\N {\bf Theorem 4.1} {\it Let $\sigma \in S_p$ and $\tau \in S_q$ be two permutations. The product $\sigma *\tau$ is
the sum of all permutations $\omega \in S_{p+q}$ verifying 

\N $\sigma \times \tau \leq \omega \leq \xi _{p,q}\cdot (\sigma
\times \tau )$, in other words:
$$\sigma *\tau =\sum _{\sigma /\tau \leq \omega \leq \sigma \backslash \tau } \omega.
$$}

\N {\it Proof.} Lemma 1.5 implies that
$$\sigma \times \tau \leq \delta \cdot (\sigma \times \tau )\leq \xi _{p,q}\cdot (\sigma \times \tau ),\ {\rm for\ any}\ \delta 
\in Sh(p,q).$$

Suppose that $\omega \in S_{p+q}$ verifies $\sigma \times \tau \leq \omega \leq \xi _{p,q}\cdot (\sigma \times \tau )$.
Let $\omega _1\in S_{p+q}$ be such that $\omega = \omega _1\cdot (\sigma \times \tau )$. It is obvious that $1_{p+q}\leq 
\omega _1$. 

Since $\omega \leq \xi _{p,q}\cdot (\sigma \times \tau )$, the definition of the weak ordering implies that there exists 
$\epsilon \in S_{p+q}$ such that $\xi _{p,q}\cdot (\sigma \times \tau )= \epsilon \cdot \omega _1\cdot (\sigma \times \tau )$, 
with $l(\xi _{p,q})=l(\epsilon )+l(\omega _1)$. It implies $\omega _1\in Sh(p,q)$, which ends the proof.\hfill  $\square $ 
\M

\N {\bf Definition 4.2} For $p,q\geq 0$, the subsets $Sh^1(p,q)$ and $Sh^2(p,q)$ of $Sh(p,q)$ are defined 
by:
$$Sh^1(p,q):=\lbrace \omega \in Sh(p,q) \mid  \omega (p+q)=p+q\rbrace ,\ {\rm and}
$$
$$Sh^2(p,q):=\lbrace \omega \in Sh(p,q) \mid  \omega (p)=p+q\rbrace .
$$

\N {\bf Remark 4.3} The set $Sh(p,q)$ is the disjoint union of $Sh^1(p,q)$ and $Sh^2(p,q)$. Moreover, one has that 
$$Sh^1(p,q)=\lbrace \omega\times 1_1 \mid  \omega \in Sh(p,q-1)\rbrace = Sh(p,q-1)\times 1_1;\ {\rm and}$$
$$Sh^2(p,q)=\lbrace (\omega \times 1_1)\cdot (1_{p-1}\vee 1_q)\mid \omega \in Sh(p-1,q)\rbrace =(Sh(p-1,q)\times
1_1)\cdot  (1_{p-1}\vee 1_q).$$
\M

\N {\bf Definition 4.4} The products $\prec $ and $\succ $ in $\overline {{\bf Q}[S_{\infty }]}$ are defined as
follows:

$$\sigma \prec \tau :=\sum _{\omega \in Sh^2(p,q)}\omega \cdot (\sigma \times \tau ),\ {\rm and}$$
$$\sigma \succ \tau :=\sum _{\omega \in Sh^1(p,q)}\omega \cdot (\sigma \times \tau ),$$
for $\sigma \in S_p$ and $\tau \in S_q$.
\M

From Remark 4.3 one gets that the associative product $*$ of {\bf Q}$[S_{\infty }]$ verifies 
$$\sigma *\tau =\sigma \prec \tau + \sigma \succ \tau ,\ {\rm for\ }\sigma ,\tau \in \overline {{\bf Q}[S_{\infty
}]}.$$

\N {\bf Proposition 4.5} {\it The operations $\prec$ and $\succ$ satisfy the relations 
$$
\eqalign{
\hbox{(i)}\quad &(a\prec b)\prec c=a\prec(b\prec c)+a\prec(b\succ c),\cr
\hbox{(ii)}\quad& a\succ (b\prec c)=(a\succ b)\prec c,\cr
\hbox{(iii)}\quad& a\succ (b\succ c)=(a\prec b)\succ c+(a\succ b)\succ c,\cr}
$$
for any $a,b,c \in \overline {{\bf Q}[S_{\infty }]}$. Hence  $\overline {{\bf Q}[S_{\infty }]}$ is a dendriform
algebra
 (as defined in [5]).}
\S
\N {\it Proof.} This is a consequence of the associativity property of the shuffle together with an inspection about
the first element of the image of the permutations. \hfill $\square$
\M

The products $\prec $ and $\succ $ may also be described in terms of the order $\leq $ as
follows.
\M

\N {\bf Proposition 4.6} {\it For any $\sigma \in S_p$ and any $\tau \in S_q$, one has:
$$\sigma \prec \tau =\sum _{(1_{p-1}\vee 1_q)\cdot (\sigma \times \tau)\leq \omega \leq \sigma \backslash \tau }\omega,\ 
$$
and
$$\sigma \succ \tau =\sum _{\sigma /\tau \leq \omega \leq (\xi _{p,q-1}\times 1_1)\cdot (\sigma \times \tau)}.$$}

\N {\it Proof.} Lemma 1.3 and Remark 4.3 imply that 
$$Sh^1(p,q)=\lbrace \omega \in S_{p+q}\mid  \omega \leq \xi _{p,q-1}\times 1_1\rbrace ,\ {\rm and}$$
$$Sh^2(p,q)=\lbrace \omega \in S_{p+q}\mid 1_{p-1}\vee 1_q\leq \omega \leq (\xi _{p-1,q}\times 1_1)\cdot (1_{p-1}
\vee 1_q)\rbrace .$$
The result follows immediately from Lemma 1.5. \hfill $\square $
\BB

\N {\bf 5. The graded algebra of planar binary trees Q$[Y_{\infty }]$.}
\M

The graded vector space {\bf Q}$[Y_{\infty }]:=\oplus _{n\geq 0}${\bf Q}$[Y_n]$ is a graded associative algebra 
for the product $*$ defined recursively as follows :
\S

\N - $t*\vert = \vert * t := T$, for all $t\in Y_n,\ n\geq 1$,

\N - if $t=t^l\vee t^r$ and $w=w^l\vee w^r$, then
$$t*w:=(t*w^l)\vee w^r + t^l\vee (t^r*w).
$$
\M

Moreover, the map $\psi^{*} :${\bf Q}$[Y_{\infty }]\rightarrow ${\bf Q}$[S_{\infty }]$, defined by
$$
\psi_n^*(t) := \sum _{\psi_n(\sigma) =t} \sigma ,
$$ is an algebra homomorphism (cf. [6]).
\M

\N {\bf Theorem 5.1} {\it If $t$ and $w$ are two planar binary trees, then the product $t*w$ verifies
$$t*w=\sum _{t/ w\leq u\leq t\backslash w}u.
$$}
\S

\N {\it Proof.} Since the ordering  $\leq $ on $Y_n$ is induced by the weak ordering  of $S_n$,
the result is a straightforward consequence of Proposition 2.8 and Theorem 4.1. \hfill $\square $
\M

As in the case of the algebra {\bf Q}$[S_{\infty }]$, we may describe on $\overline {{\rm {\bf Q}}[Y_{\infty }]}
:=\oplus _{n\geq 1}${\bf Q}$[Y_n]$ two products $\prec $ and $\succ $, such that:
$$t*w=t\prec w+t\succ w,\ {\rm for\ any\ }t,w\in \overline {{\rm {\bf Q}}[Y_{\infty }]}.
$$
\M

\N {\bf Definition 5.2\ }ÊLet $t\in Y_p$ and $w\in Y_q$. The elements $t\prec w$ and $t\succ w$ in $\overline {{\rm
{\bf Q}}[Y_{\infty }]}$  are given  by:
$$t\prec w:=t^l\vee (t^r*w),\ {\rm for\ }t=t^l\vee t^r,$$
$$t\succ w:=(t*w^l)\vee w^r,Ö {\rm for\ }w=w^l\vee w^r.$$
\M

The space $\overline {{\rm {\bf Q}}[Y_{\infty }]}$, equipped with the products $\prec $ and $\succ $ is a dendriform algebra 
(cf. [5], [6]). We prove now that $\overline {\psi }:\overline {{\rm {\bf Q}}[S_{\infty }]}\rightarrow 
\overline {{\rm {\bf Q}}[Y_{\infty }]}$ preserves $\prec $ and $\succ $.
\M 

\N {\bf Proposition 5.3} {\it The $K$-linear map $\overline {\psi ^*}:\overline {{\rm {\bf Q}}[Y_{\infty }]}\rightarrow 
\overline {{\rm {\bf Q}}[S_{\infty }]}$ is a dendriform algebra homomorphism.}
\S

\N {\it Proof.} We prove that $\psi^*(t\succ w)=\psi^*(t)\succ \psi^*(w)$, for any trees $t $ and $w$. The proof
that $\psi ^*$ preserves the product $\prec $ is analogous.

Recall that the associativity of the shuffle product is equivalent to the following equality:
$$Sh(p,q+r)\cdot (1_p\times Sh(q,r))=Sh(p+q,r)\cdot (Sh(p,q)\times 1_r).
$$
Let $t\in Y_p$, and $w=w^l\vee w^r\in Y_{q+r}$ with $w^r\in Y_q$ and $w^l\in Y_r$. Recall from [5] that 
the right product is given by
$t\succ w = (t*w^l)\vee w^r$. So,
$$\psi^*(t\succ w) = \psi^*((t*w^l)\vee w^r) = \sum_{\gamma \in Sh(p+q-1,r)} (\gamma\times 1_1) \cdot 
\big( \psi^* (t*w^l)\vee \psi^* ( w^r)\big).
$$
Since 
$$\psi^*(t* w^l) = \psi^*(t)*\psi^*( w^l)= \sum_{\delta \in Sh(p,q-1)} \delta \cdot 
\big( \psi^* (t)\times \psi^* ( w^l)\big),
$$
one has by (5.3.1):
$$\displaylines{
\psi^*(t\succ w) \hfill \cr
 \hfill= \sum_{\gamma \in Sh(p+q-1,r)} \sum_{\delta \in Sh(p,q-1)}
(\gamma\times 1_1)\cdot  (\delta\times 1_{r+1})\cdot \big( (\psi^*(t)\times \psi^* ( w^l))\vee  \psi^*(
w^r)\big)\cr
\hfill=  \sum_{\omega \in Sh(p,q+r-1)} \sum_{\epsilon \in Sh(q-1,r)}
(\omega\times 1_1)\cdot  (1_p\times \epsilon\times 1_1)\cdot \big( (\psi^*(t)\times \psi^* (w^l))\vee  \psi^*
( w^r)\big).\cr
}$$
Since 
$$\displaylines{
 (\psi^*(t)\times \psi^*(w^l))\vee  \psi^*(w^r))= \big(\psi^*(t)\times \psi^*(w^l)\times \psi^*(w^r)
\times 1_1\big) \cdot s_{p+q+r-1}\cdots s_{p+q} \cr
\hfill =  \psi^*(t)\times \big(\psi^*(w^l))\vee  \psi^*(w^r)\big),\cr
}$$
we get:
$$\psi^*(t\succ w) =  \sum_{\omega \in Sh(p,q+r-1)} (\omega\times 1_1)\cdot \big( \psi^*(t)\times 
 \sum_{\epsilon \in Sh(q-1,r)}  (\epsilon\times 1_1)\cdot \big( \psi^*(w^l)\vee  \psi^*(w^r)\big)
$$
$$\hfill =  \sum_{\omega \in Sh(p,q+r-1)} (\omega\times 1_1)\cdot \big( \psi^*(t)\times \psi^*(w))= \psi^*(t)\succ
\psi^*(w).
$$
{ }\hfill $\square$
\BB

\N {\bf 6. The graded algebra of the cube vertices {\bf Q}$[Q_{\infty }]$.}
\M

The image of $(\phi_n \circ \psi_n )^*$ in {\bf Q}$[S_n]$ is the so-called {\it Solomon algebra}, that   we 
denote  by  $Sol_n$.
The direct sum {\bf Q}$[Q_{\infty }]:= \bigoplus_{n\geq 0}{\bf Q}[Q_n]$ is a graded subalgebra of {\bf Q}$[Y_{\infty
}]$ and so of {\bf Q}$[S_{\infty }]$. Since, by section 3, $\psi_n$ is compatible with the orders and with the `over
' and `under' operations, the same  kind of arguments as in section 5 implies the following result:
\M 

\N {\bf Theorem 6.1.} {\it For any $\epsilon\in Q_p$ and any $\delta\in Q_q$, the product $*$ verifies:
$$ \epsilon * \delta= \sum_{\epsilon / \delta\leq \alpha \leq  \epsilon \backslash \delta} \alpha \  = \epsilon /
\delta + \epsilon
\backslash \delta.
$$}
Recall from section 3 that:
$$\displaylines{
\epsilon / \delta := (\epsilon _1,\dots ,\epsilon _{p-1},-1,\delta _1,\dots ,\delta _{q-1})\cr
\epsilon \backslash \delta := (\epsilon _1,\dots ,\epsilon _{p-1},+1,\delta _1,\dots ,\delta _{q-1}).\cr
}$$
Since there is obviously no element in between 
$\epsilon / \delta$ and $\epsilon \backslash \delta$ the formula for the product on the generators takes the
  form $\epsilon *\delta = \epsilon / \delta + \epsilon \backslash \delta$, that is
$$\displaylines{
 (\epsilon _1,\dots ,\epsilon _{p-1})*(\delta _1,\dots ,\delta _{q-1})= \hfill\cr
\hfill (\epsilon _1,\dots ,\epsilon _{p-1},+1,\delta _1,\dots ,\delta _{q-1})+
(\epsilon _1,\dots ,\epsilon _{p-1},-1,\delta _1,\dots ,\delta _{q-1}).\cr
}$$
Hence we recover exactly formula 4.6 of [6, p. 307].
\BB 

\N {\bf Appendix. The weak Bruhat order  on a Coxeter group.}
\M

Let $(W,S)$ be a finite Coxeter system (cf. [1]). So $W$ is a finite group generated by the set $S$, with relations
of the form
$$(s\cdot s\rq )^{m(s,s\rq )}=1,\ {\rm for}\ s,s\rq \in S,$$
for certain positive integers $m(s,s\rq )$, with $m(s,s)=1$ for all $s\in S$.

For any element $w\in W$  the length  $l(w)$ is the number of factors in a minimal expression of $w$ in terms of elements in $S$.
There  exists a unique element of maximal length in $W$, denoted $w^0$.

Given a subset $J\subseteq S$, the {\it standard parabolic subgroup} $W_J$ is the subgroup of $W$ generated by $J$. Clearly, 
the pair $(W_J,J)$ is a finite Coxeter system too.

\N {\bf Definition A.1} Let $(W,S)$ be a finite Coxeter system and let $J$ be a subset of $S$. The set $X_J$ of elements
 of $W$ that have no descent at $J$  is defined as
$$X_J:=\lbrace w\in W \mid  l(w\cdot s)>l(w),\ {\rm for\ all}\ s\in J\rbrace .
$$

The following result is due to L. Solomon.
\M

\N {\bf Proposition A.2} ([7] p. 258) {\it Let $(W,S)$ be a finite Coxeter system, and let $J$ be a subset
 of $S$. Every
element of $W$  can be written uniquely as $w=x\cdot y$, where $x\in X_J$ and $y\in W_J$.
 If $x\in X_J$ and $y\in W_J$,
then $l(x\cdot y)=l(x)+l(y)$.}
\M

\N {\bf Definition A.3} Let $(W,S)$ be a finite Coxeter system, the {\it weak Bruhat order} on $W$ is defined by:
$$x\leq x\rq \ {\rm if}\ x=y\cdot x\rq,\ {\rm with}\ l(x)=l(y)+l(x\rq ).$$
\M

The group $W$ equipped with the weak ordering is a finite poset with minimal element $1_W$, and maximal element $w^0$.

Given a subset $J\subseteq S$, Solomon\rq s result implies that there exist unique elements $x_J^0\in X_J$ and 
$w_J^0\in W_J$ such that $w^0=x_J^0\cdot w_J^0$. It is easy to check that $w_J^0$ is the maximal element of $(W_J,J)$, and that 
$x_J^0$ is the longest element of $X_J$.
\M

\N {\bf Corollary A.4} {\it Let $(W,S)$ be a finite Coxeter system and let $J\subseteq S$, then $X_J$ is the subset of
$W$ characterized as follows:
$$X_J=\lbrace w\in W \mid  w\leq x_J^0\rbrace .
$$}

\N {\it Proof.} For any  $w\in X_J$ one has  $l(w\cdot w_J^0)=l(w)+l(w_J^0)$ and $w\cdot w_J^0\leq w^0$. There exists 
$y\in W$ such that $w^0=y\cdot w\cdot w_J^0$, with $l(w^0)=l(y)+l(w)+l(w_J^0)$, which implies $x_J^0=y\cdot w$, with 
$l(x_J^0)=l(y)+l(w)$. So one has $w\leq x_J^0$.

Conversely, if $w\leq x_J^0$, then there exists $y\in W$ such that $x_J^0=y\cdot w$, with $l(x_J^0)=l(y)+l(w)$. If $w\notin X_J$ 
there exists $s\in J$ such that $l(w\cdot s)<l(w)$, then $l(x_J^0\cdot s)=l(y\cdot w\cdot s)<l(y\cdot w)=l(x_J^0)$. 
But $x_J^0\in X_J$, so $w\in X_J$.\hfill  $\square $
\BB

\centerline {\bf References}
\M
\N [1] N. Bourbaki, {\it Groupes et alg\`ebres de Lie, Ch. 4-6}, Hermann Paris {\oldstyle 1968}; Masson, Paris,
{\oldstyle 1981}.
\S
\N [2] Chr. Brouder and A. Frabetti,   Renormalization of QED with trees,  European Physical Journal C, 
{\tt [hep-th/0003202]}.
\S
\N [3] I.M. Gelfand,  D. Krob, A. Lascoux, B. Leclerc,
 V. Retakh,  J.-Y. Thibon,   
Noncommutative symmetric functions, 
Adv. Math. 112 ({\oldstyle 1995}), no. 2, 218--348.
\S
\N [4] C. Malvenuto,  and Chr. Reutenauer,    
Duality between quasi-symmetric functions and the Solomon descent algebra,
J. Algebra 177 ({\oldstyle 1995}), no. 3, 967--982.
\S
\N [5] J.-L. Loday, Dialgebras, preprint IRMA, Strasbourg,  {\oldstyle 1999}. {\tt [math.QA/0102053]}.
\S
\N [6] J.-L. Loday,  and M. O. Ronco,  Hopf algebra of the planar binary trees,  Adv. Math. 139 ({\oldstyle
1998}), no. 2, 293--309.
\S
\N [7] L. Solomon,    
A Mackey formula in the group ring of a Coxeter group, 
J. Algebra 41 ({\oldstyle 1976}), no. 2, 255--264. 
\BB
\N JLL : Institut de Recherche Math\'ematique Avanc\'ee,

    CNRS et Universit\'e Louis Pasteur

    7 rue R. Descartes,

    67084 Strasbourg Cedex, France

    E-mail : loday@math.u-strasbg.fr
\BB

\N MOR : Departamento de Matem\'atica

Ciclo B\'asico Com\'un

Universidad de Buenos Aires

Pab. 3 Ciudad Universitaria Nu\~nez 

(1428) Buenos-Aires, Argentina

E-mail : mronco@mate.dm.uba.ar
\BB

\end